\documentclass[12pt]{article}
\usepackage[T2A]{fontenc}
\usepackage[utf8]{inputenc}
\usepackage[russian,english]{babel}
\usepackage{amsmath,amssymb,euscript,amsthm,amsfonts,mathrsfs,amscd,mathtext}
\usepackage{color}
\usepackage{cite,enumerate,float,indentfirst}
\usepackage{graphicx}
\usepackage{setspace}

\theoremstyle{plain}

\newtheorem{lem}{Lemma}%[section]

\theoremstyle{Definition}

\newtheorem{definition}{Definition}

\theoremstyle{definition}
\newtheorem{exmp}{Example}%[section]
\theoremstyle{Remark}
\newtheorem{rem}{Remark}%[section]
\newcommand{\dd}{\ensuremath{\displaystyle}}

\newcommand{\infl}{\ensuremath{\inf\limits}}

\newcommand {\intl}{\ensuremath{\int\limits}}
\newcommand{\suml}{\ensuremath{\sum\limits}}
\newcommand{\PP}{\ensuremath{\mathbf{P}}}

\newcommand{\EEE}{\ensuremath{\mathscr{E}}}
\newcommand{\PPP}{\ensuremath{\mathscr{P}}}

\newcommand{\GGG}{\ensuremath{\mathscr{G}}}
\newcommand{\FFF}{\ensuremath{\mathscr{F}}}

\newcommand{\SSS}{\ensuremath{\mathscr{S}}}
\newcommand{\UU}{\ensuremath{\mathscr{U}}}

\newcommand{\EE}{\mathbb E}

\newcommand{\bD}{\stackrel{\mathscr{D}}{=}}
\newcommand{\1}{\ensuremath{\mathbf{1}}}

\newcommand{\ud}{\,\mathrm{d}}
\newcommand*{\TR}{\hfill\ensuremath{\triangleright}}
\newcommand{\bd}{\stackrel{\text{\rm def}}{=\!\!\!=}}

\title{On coupling epoch for renewal times}
\author{Galina A. Zverkina\footnote{ Institute of Control Sciences V.A. Trapeznikov Academy of Sciences; the author is supported by RFBR, project No 20-01-00575A}}
\begin{document}

\maketitle

\section{Introduction}
It is well-known that the queuing theory and related theories (for example, reliability theory or queuing network theory) are based on the renewal theory (see, e.g., \cite{GK}).
Thus, for the study of the convergence rate of the distribution of some queuing system (or reliability system, or network) it is necessary to be able to estimate the rate of convergence of distributions of the ``{\it states}'' of renewal processes.
Indeed, most of the studied queueing processes are regenerative. 
Therefore, the distribution of a queueing process is determined by the elapsed time since the last regeneration time. 

The regenerative process has an embedded renewal process.
For regenerative processes, the elapsed time since the last regeneration time is the backward renewal time of the corresponding renewal process.

Recall that the renewal process $R _t\bd\dd\dd \suml _{i=1} ^\infty \1\left \{ \dd \suml _{s=1} ^i \xi _k\leqslant t\right \} $, where
$\left \{\xi _1, \xi _2, ...\right \} $ are i.i.d. positive random variables.
$R _t$ is a counting process which changes its states at times $t _k=S _k\bd\dd\dd \suml _{j=1} ^s \xi _j$.
The times $t _k$ are renewal times.

\begin{figure}[h]
\centering
\begin{picture}(300,45)
\put(0,20){\vector(1,0){300}}
\put(0,18){\line(0,1){4}}
\thicklines
{\qbezier(0,20)(30,40)(60,20)}%zeta0
{\qbezier(60,20)(80,40)(100,20)}
\qbezier(100,20)(150,40)(200,20)%zeta3
\qbezier(200,20)(240,40)(280,20)%zeta4
%\qbezier(240,20)(260,0)(280,20)%zeta4
%\qbezier(280,20)(320,30)(330,31)%zeta4
\put(280,30){\ldots \ldots \ldots}
\put(58,10){$t _1$}
\put(98,10){$t _2$}
\put(198,10){$t _3$}
\put(278,10){$t _4$}
\put(27,40){$\xi _1$}
\put(77,40){$\xi _2\bD \xi _1$}
\put(147,40){$\xi _3\bD \xi _1$}
\put(237,40){$\xi _4\bD \xi _1$}
%\multiput(0,5)(3,0){11}{\ww \line(0,1){40}}
\put(-5,10){$t _0=0$}
\put(148,10){$t$}
\put(195,15){\line(1,1){10}}
\put(195,25){\line(1,-1){10}}
%
%\rr
\qbezier(150,20)(175,0)(200,20)
\qbezier(150,20)(175,1)(200,20)
\qbezier(150,20)(175,2)(200,20)
\put(170,0){$W _t$}
%
%
%\bb
\qbezier(100,20)(125,0)(150,20)
\qbezier(100,20)(125,1)(150,20)
\qbezier(100,20)(125,2)(150,20)
\put(120,0){$B _t$}
%
%\multiput(102,2)(3,0){25}{объект}
%\bb
%\qbezier(243,20)(271,0)(300,20)
%\put(270,0){$w' _t$}
%\put(-5,10){$-b' _0$}
%\thinlines
\end{picture}
\caption{$B _t$ is a backward renewal time, and $W _t$ is a forward renewal time at the fixed time $t$; $B _t\bd t-\dd\dd \suml _{1} ^{N _{t}}\xi _i;$ $W _t\bd t-S _{N _t}$.}
\label{fig1}%\mbox{\bf ~~~ fig1}
\end{figure} 

Fig.\ref{fig1} shows the {\it backward renewal time (or overshoot)} $B _t$, and the {\it forward renewal time (or undershot)} $W _t$ at the {\bf fixed} time $t$.

If the distribution $\PPP(B_t)$ of the process $B_t$ converges weakly to the stationary distribution $\tilde\PPP_B$ as $t\to \infty$ ($\PPP(B_t) \Rightarrow \tilde\PPP_B$), then $\|\PPP(B_t)-\tilde\PPP_B\|_{TV}\to 0$.

The value of $B_t$ defines the distribution $\PPP(X_t)$ of corresponding regenerative process $X_t$ at the time $t$. 
Hence, $\PPP(X_t) \Rightarrow \tilde\PPP_X$, and $\|\PPP(X_t)-\tilde\PPP_X\|_{TV} \leqslant \|\PPP(B_t)-\tilde\PPP_B\|_{TV} $.
Here $\tilde\PPP_Y$ denotes the stationary distribution of the process $Y_t$.

Recall the 

\begin{definition}
The total variation distance between two probability measures $\PPP_1$ and $\PPP_2$ on a sigma-algebra $\FFF$ of subsets of the sample space ${ \Omega }$ is defined via
${\displaystyle \delta (\PPP_1,\PPP_2)=\|\PPP_1-\PPP_2\|_{TV}=\sup _{A\in {\FFF}}\left|\PPP_1(A)-\PPP_2(A)\right|.}
$
\TR
\end{definition}

In the papers \cite{zvconver,zvlor}, the method of the construction of upper bounds for $\|\PPP(B_t)-\tilde\PPP_B\|_{PV}$ was given.
However, in these papers, the general interest is the behaviour of regenerative processes.
Thus, the calculation algorithm for the upper bounds is not complete.

The aim of the presented paper is a detailed description of this algorithm.

\section{On coupling epoch for renewal times of two renewal processes $R_t$ and $R_t'$}

So, we have the regenerative process $X_t$, corresponding embedded renewal process $R_t$, and $R_t$ generates the Markov process $B_t$.
The renewal periods of $R_t$ are i.i.d. r.v.'s $\xi_i$ with d.f. $F(s)$; for simplicity, we suppose here that d.f. $F(s)$ is absolutely continuous.

In this paper, we suppose that:
\begin{itemize}
 \item There exists at least two finite moments of r.v. $\xi_i$.
 \item For all $s \geqslant 0$, probability density function (p.d.f.) $f(s)\bd F'(s)$ of r.v. $\xi_i$ is positive almost surely (a.s.).
\end{itemize}

It is well-known that $\PP\{B_{t+\Delta} \geqslant a+ \Delta | B_t>a\in \}= 1-\lambda(a)\Delta + o(\Delta)$, where $\lambda(s)\bd \dd\dd \frac{F'(s)}{1-F(s)}$ -- the {\it intensity} of renewals.

The basic schema and some conditions of the construction of the upper bounds for the coupling epoch of the renewal times for two renewal processes with different initial values of backward renewal time we're given in papers \cite{zvconver,zvlor}.

Let $R_t$ and $R_t'$ be two renewal processes (with the same distributions of renewal periods) that started {\it before} the time $t_0=0$.
Thus, at the time $t_0=0$ corresponding backward renewal times $B_0=b$, $B_0'=b'$ $\in [0,\infty)$.
Denote by $t_i$, $t_i'$ consecutive recovery times of processes 
Put $T\bd\{t_1, t_1'\}$, where $t_i$, $t_i'$ are consecutive renewal times of the processes $R_t$, $R_t'$ accordingly.
The distributions of $t_1$ and $t_1'$ are the distributions of the forward renewal times of the processes $R_t$ and $R_t'$ at the time $t=0$.

Thus, $\PP\{t_1>s\}=\dd\dd \frac{1-F(s+b)}{1-F(b)}$, and $\PP\{t_1'>s\}=\dd\dd \frac{1-F(s+b')}{1-F(b')}$.
So, $\EE\,\varphi(\max\{t_1,t_1'\}) \leqslant \EE\,\varphi(t_1+t_1') $ for any increasing function $\varphi(s)$.

For simplicity suppose that $T=t_1$.

 In the above conditions, at the times $t_1$, $t_2$, \ldots, $t_i$, \ldots, $\PP\{B_t' \leqslant \Theta\} \geqslant p_0=p_0(\Theta)=1-\dd\dd \frac{R(\xi)}{\Theta}$, where $R(\xi)\bd\dd\dd \frac{\EE\,\xi^2}{\EE\,\xi}$, $\Theta$ is some number from $( R(\xi),+\infty)$. It was proved in \cite{zvconver,zvlor}.
Hereinafter $\xi$ is an r.v. with d.f. $F(s)$.
The calculation of the bound $p_0$ is based on Lorden's inequality.

Denote the event $\{B'_{t_i}<\Theta\}$ by $\SSS_i$; $\PP(\SSS_i) \geqslant p_0$.

Thus, we can apply the {\it Basic Coupling Lemma} (BCL) at times $t_i$:
\begin{lem}
Let $f_i(s)$ be the distribution density of r.v. $\theta_i$ ($i=1,2$).
Let $\varphi(s)\in[0,\min(f_1(s),f_2(s))]$, and
$$\dd \intl_{-\infty}^ \infty \varphi(s) \ud s=\varkappa>0.$$

Then on some probability space there exists two random variables $\vartheta_i$ such that $\vartheta_i\bD \theta_i$, and $\PP\{\vartheta_1=\vartheta_2\}\geqslant \varkappa$. \hfill \ensuremath{\triangleright}
\end{lem}
The proof of BCL is very easy, but further, we will use elements of this proof and the corollary from this Lemma. 
\begin{proof}
Denote $\varphi(s)\stackrel {\text{\rm def}}{=\!\!\!=} \min\left(f_1(s), 2_2(s)\right)$; $\varkappa(\varphi) \bd \dd\int\limits_0^\infty \varphi(s)\,\mathrm{d}\, s =\Phi(+\infty)>0$ where $\Phi(s)\stackrel {\text{\rm def}}{=\!\!\!=} \dd\int\limits_0^s \varphi(u)\,\mathrm{d}\, u$.
Put $F_i(s)\stackrel {\text{\rm def}}{=\!\!\!=} \dd\int\limits_0^s f_i(u)\,\mathrm{d}\, u$.

Denote
$\Psi_i(s)\stackrel {\text{\rm def}}{=\!\!\!=} F_i(s)-\Phi(s)$, $\Psi_i(+\infty)=1-\varkappa(\varphi)$.

Put for i.i.d. r.v. $\UU$, $\UU'$, $\UU''$ uniformly distributed on $[0,1)$
$$
\vartheta_i(\UU,\UU',\UU'')\stackrel {\text{\rm def}}{=\!\!\!=} \mathbf{1} (\UU<\varkappa(\varphi))\Phi^{-1}(\varkappa(\varphi) \UU')+\mathbf{1} (\UU\geqslant \varkappa(\varphi)) \Psi_i^{-1}((1-\varkappa(\varphi)) \UU'').$$

It is easy to see that 
 $\PP\Big\{\vartheta_i(\UU,\UU',\UU'')\leqslant s\Big\}=F_i(s)$, and\\
$\PP\left\{\vartheta_1(\UU,\UU',\UU'')=\vartheta_2(\UU,\UU',\UU'')\right\}=\varkappa(\varphi)$.
\end{proof}
\begin{rem}
If $\varphi(s)\equiv \varphi_0(s)\bd \min(f_1(s),f_2(s))$, then \linebreak $\PP\left\{\vartheta_1(\UU,\UU',\UU'')=\vartheta_2(\UU,\UU',\UU'')\right\}=\varkappa_0$. \TR
\end{rem}
The value $\varkappa_0$ is called a {\it common part of distributions of $\theta_1$, $\theta_2$. }
\begin{rem}
For any non-negative function $\varphi(s)$, 
\begin{equation}\label{Efi}%{\bb Efi\;\;\;}
 \EE\,\varphi(\vartheta_1) \leqslant \varkappa\,\EE \,\varphi(\vartheta_i); \qquad \EE\,\varphi(\vartheta_2) \leqslant (1-\varkappa)\,\EE \,\varphi(\vartheta_i). 
\end{equation}
{\TR}
\end{rem}

So, at the time $t_i$ in the case when the event $\SSS_i$ happened, $B'_{t_i}=\theta\in[0,\Theta)$, and the distribution of the forward renewal time $W_{t_i}'$ is $F^W_\theta(s)=1-\dd\dd \frac{1-F(s+\theta)}{1-F(\theta)}$; denote $f^W_\theta(s)\bd \dd \frac{\ud}{\ud s}F^W_\theta(s)$.

So, the common part of the distributions of $\xi_i$ and $W'_{t_i}$ is \begin{multline*}
 \varkappa(\theta)\bd \dd \intl_{-\infty}^ \infty \min(f(s),f^W_\theta(s)(s))\ud s  \geqslant
 \\
  \geqslant \varkappa_\Theta\bd \infl{\theta\in[0,\Theta)}\dd \intl_{-\infty}^ \infty \min(f(s),f^W_\theta(s)(s))\ud s >0. 
\end{multline*}
Thus, at the time $t_i$, we can prolong the processes $R_t$ and $R_t'$ by such a way, that the time $t_{i+1}=t_i+\xi_i$ is the coupling epoch $\tau(b,b')$ (the time of the coincidence of the processes $B_t$ and $B_t'$) with probability greater then $\pi\bd p_0(\Theta)\varkappa(\Theta)$.
Moreover, the marginal distributions of the created and original processes remain unchanged.

Denote $\EEE_i\bd \{ \SSS_i\cap \{\tau(b,b')=t_{i+1}\}\}$, and $\GGG_i\bd \EEE_i\cap \bigcap\limits_{j=1}^{i-1}\overline{\EEE_j}$, $j\in \mathbb{N}$; here we put $\bigcap\limits_{j=1}^{0}\overline{\EEE_j}\bd \varnothing$.
In these denotations, $\tau(b,b')$ is less than the {\it conditional} geometrical sum of independent random variables:
$\tau \leqslant T+\dd \suml_{i=1}^{\nu} \{\xi_i|\GGG_\nu \} $, and $\PP\{\nu>n\} \leqslant (1-\pi)^n$ because $\PP(\overline{\EEE_i}) \leqslant q\bd (1-\pi)$, $\PP({\EEE_i}) \leqslant 1$.

\section{Polynomial bounds for the coupling epoch}
Let in addition for the conditions of the previous part, $k \geqslant 2$ finite moments of $\xi_i$ exist.

For estimate $\tau(b,b')^\ell$ for all $\ell\in [1,k]$ we will use the Jensen's inequality:

For a real convex function $\varphi$ , numbers $ x_{1},x_{2},\ldots ,x_{n}$ in its domain, and positive weights $a_{i}$, Jensen's inequality can be stated as:

$$
\varphi \left({\dd \frac {\sum a_{i}x_{i}}{\sum a_{i}}}\right)\leq {\dd \frac {\sum a_{i}\varphi (x_{i})}{\sum a_{i}}}
\mbox{\; or \;}
\varphi \left({\dd \frac {\sum x_{i}}{n}}\right)\leq {\dd \frac {\sum \varphi (x_{i})}{n}}.
$$

If there exists $\EE\,\xi_i^{\ell}<\infty$, we estimate:
\begin{multline*}
 \EE\left(T_1+\dd \suml_{i=1}^{\nu} \xi_i\right)^\ell=\EE\left(\dd \frac{(\nu+1)\times\left( T+\dd \suml_{i=1}^{\nu} \{\xi_i|{\GGG_\nu} \}\right)}{\nu+1}\right)^\ell \leqslant
 \EE\left((\nu+1)^\ell\dd \frac{\left(T_1\right)^\ell+\dd \suml_{i=1}^{\nu} \{\xi_i^\ell|{\GGG_\nu} \}}{\nu+1}\right) \leqslant
 \\
  \leqslant
 \EE\left((\nu+1)^{\ell-1}\right)\EE\,T_1^\ell +\dd \suml_{i=1}^\infty \left(\left(\dd \suml_{j=1}^{i-1}[\EE(\xi_j^\ell|\overline{\EEE_i})\PP(\overline{\EEE_j})]q^{i-2}\right)+\EE[(\xi_i^\ell|\EEE_i)\PP(\EEE_i)]q^{i-1} \right) \leqslant
 \\
  \leqslant
 \EE\left((\nu+1)^{\ell-1}\right)\EE\,T_1^\ell +\dd \suml_{i=1}^\infty \left(\dd \suml_{j=1}^{i-1}[\EE(\xi_j^\ell)]q^{i-2}+\EE(\xi_i^\ell)q^{i-1} \right) \leqslant
 \\
  \leqslant
 \EE\,(t_1+t_1')^\ell\times \left(\dd \suml_{i=0}^\infty (i+1)^{\ell-1}\times q^i\right) + \EE (\xi^\ell)\times\left(\left(\dd \suml_{i=1}^\infty i\times q^{i-1}\right)+\dd \suml_{i=0}^\infty q^i\right)=
 \\
 =\EE\,(t_1+t_1')^\ell\times S_\ell+ \EE(\xi^\ell)\times \left(\dd \frac{1}{1-x^2}+\dd \frac{1}{1-x}\right)=\mbox{\bf Poly}(\tau(b,b'),\ell),
\end{multline*}
where $S_\ell$ can be calculated by the formula $\dd \suml_{k=0}^\infty k^\ell x^k = \left(x \dd \frac{\ud}{\ud x}\right)^\ell\dd \frac{1}{1-x}$ for natural $n$.

For non-integer $\ell$, the upper bounds for $S_\ell$ can be calculated numerically.

In formulae above, $\dd \suml_1^0 (\cdot)=\dd \suml_1^{-1} (\cdot)\bd 0$.

\section{Exponential bounds for the coupling epoch}
In addition for the conditions of previous part let for some $\alpha>0$, $\EE\,\exp(\alpha \xi_i)=A<\infty$.

In this case we know that $\EE\,\exp(\alpha (\xi_i|\EEE_i))\PP(\EEE_i) \leqslant \EE\,\exp(\alpha \xi_i)>1$. 

Thus, for finding upper bounds for the sum 
\begin{multline}\label{prod}%{\bb prod\;\;\;}
\EE\,\exp(\beta \tau(b,b'))= \EE\,\exp(\beta \dd \suml_{j=1}^\nu \{\xi_i|\GGG_\nu\})= 
\\
=\EE \,\exp(\beta\cdot T)\left(\prod\limits_{j=1}^{\nu-1}\PP(\overline{\EEE}_j) \EE\,\exp(\beta\cdot \{\xi_j|\overline{\EEE}_j\})\right)\times \PP({\EEE}_\nu) \EE\,\exp(\beta\cdot \{\xi_i|{\EEE}_\nu\}) \leqslant
\\
 \leqslant
\EE \,\exp(\beta\cdot T)\EE\,\exp(\beta\cdot \xi)\dd \suml_{i=1}^\infty q^{i-1}\prod\limits_{j=1}^i\EE\,\left( \exp(\beta \{\xi_i|\overline{\EEE}_j\})\right),
\end{multline}
here anew, $\dd \suml_1^0 (\cdot)\bd 0$, $\prod\limits_1^{0} (\cdot)\bd 1$, and $\xi$ is an r.v. with d.f. $F(s)$.

Note that $\EE\, \exp(\beta \{\xi_i|\overline{\EEE}_j\})>1$.

There is a number $\beta\in (0,\alpha)$ such that $q\EE\, \exp(\beta \{\xi_i|\overline{\EEE}_i\})<1$, and the sum (\ref{prod}) can be calculated.

But now we do not know the general way to find such a number $\beta$.
For any cases, we need to use peculiar properties of d.f. $F(s)$.

\begin{exmp}
Suppose that d.f. $F(s)$ is defined by its intensity $\lambda(s)=C+\dd\dd \frac{K}{1+s}$ \; ($C>0, K>0$).
As $\lambda(s)>C$, $\EE\,\exp(\beta \xi)<\infty$ for all $\alpha\in (0,C)$. 

$F(s)=1-\exp\left( -\dd\dd \intl_0^s \lambda(u) \ud u\right)=1-\dd\dd \frac{e^{-C s}}{(1+s)^K}$, i.e. $\xi=\min(\zeta_1,\zeta_2)$, where $\PP\{\zeta_1 \leqslant s\}=1-\exp(-Cs)$, $\PP\{\zeta_2 \leqslant s\}=1-\dd\dd \frac{1}{(1+s)^K}$.

It's p.d.f. $f(s)=F'(s)=\dd\dd \frac{e^{-Cs}(C+K+Cs)}{(1+s)^{K+1}}$.

Further, $\EE\,\xi<\min\{\dd \frac{1}{C},\dd \frac{1}{K-1}\}=\Gamma_1$, and $\EE\,\xi^2 \geqslant \EE\,\zeta_1^2+\EE\,\zeta_2^2=\dd \frac{2}{C^2}+\dd \frac{2}{(K-1)(K-2)}=\Gamma_2$; $R(\xi)>\hat R(\xi)\bd \dd\dd \frac{\Gamma_2}{\Gamma_1}$. 

Let $\Theta>\hat R(\xi)$. 
With probability greater then $p_0=1-\dd\dd \frac{\hat R(\xi)}{\Theta}$ the event $\{B_t' \leqslant \Theta\}$ happens.

The conditional distribution $\PP\{W_t' \leqslant s|B_t'=\theta\}=1-\dd\dd \frac{1-F(s+\theta)}{1-F(\theta)}$; the conditional density of $\{W_t'|B_t'=\theta\}$ is $f^W(s)=\dd\dd \frac{e^{-Cs}(C+K+C(s+\theta))(1+\theta)^K}{(1+s+\theta)^{K+1}}$.

Therefore (see the Proof of BCL above), $\varphi_0(s)\bd\min\{f(s), f^W_\theta(s)\} \geqslant \varphi(s)\bd \dd\dd \frac{Ce^{-Cs}}{(1+s+\Theta)^{K+1}}$, the value $\varkappa(\varphi)$ can be calculated numerically.

A rough upper bound of this value is
$$
\varkappa(\varphi)=\dd \intl_0^\infty \dd \frac{Ce^{-Cs}}{(1+s+\Theta)^{K+1}}\ud s < \dd \intl_0^\infty \dd \frac{Ce^{-Cs}}{(1+\Theta)^{K+1}}\ud s=\dd \frac{1}{(1+\Theta)^{K+1}}, 
$$
and $1-\varkappa(\varphi)>1-\dd\dd \frac{1}{(1+\Theta)^{K+1}}$.

Thus (see the Proof of BCL above), 
\begin{multline*}
\psi(s)=\dd \frac{\ud }{\ud s}\PP\{W_\theta \leqslant s\}-\varphi(\theta)=
\\
=\dd\dd \frac{e^{-Cs}(C+K+C(s+\theta))(1+\theta)^K}{(1+s+\theta)^{K+1}}-\dd\dd \frac{Ce^{-Cs}}{(1+s+\Theta)^{K+1}} <
\\
< 
\dd\dd \frac{e^{-Cs}(C+K+C(s+\Theta))(1+\Theta)^K}{(1+s)^{K+1}}-\dd\dd \frac{Ce^{-Cs}}{(1+s)^{K+1}}=
\\
=e^{-Cs}\dd \frac{(C+K+C(s+\Theta))(1+\Theta)^K-C}{(1+s)^{K+1}}.
\end{multline*}
So, for r.v. $\vartheta$ with d.f. $\Psi(s)\bd\dd \intl_0^\infty \dd \frac{\psi(s)}{1-\varkappa(\varphi)}\ud s$ we can estimate:
\begin{multline*}
\EE\,\exp(\beta \{W_t' \leqslant s|B_t'=\theta\}) \leqslant 
\\
 \leqslant\dd \frac{1}{1-\varkappa(\varphi)}\dd \intl_0^\infty e^{\beta s} e^{-Cs} \dd \frac{(C+K+C(s+\Theta))(1+\Theta)^K-C}{(1+s)^{K+1}} \leqslant 
\\
 \leqslant
\left(\dd \intl_0^\infty e^{-2(C-\beta)s}\ud s \dd \intl_0^\infty \left(\dd \frac{(C+K+C(s+\Theta))(1+\Theta)^K-C}{(1+s)^{K+1}}\right)^2 \right)^{\dd \frac12}
=
\\
=
\sqrt{\dd \frac{1}{2(C-\beta)}}\,Q(\Theta),
\end{multline*}
i.e. there exists $\EE\,\exp(\beta \{W_t' \leqslant s|B_t'=\theta\})$ for all $\beta\in (0,C)$. 

Finally, by analytical or numerical methods, we can find the number $\beta_0>0$ such that $\EE\,\exp(\beta \{W_t' \leqslant s|B_t'=\theta\}) \leqslant M(\beta,\Theta)<\dd\dd \frac{1} {q}$ (about $q$ see the Section 2) for all $\beta<\beta_0$, $\theta\in[0,\Theta]$.

Therefore, for these numbers $\beta$ we have: 
\begin{multline*} 
\EE\,\exp(\beta\tau(b,b'))  \leqslant \EE \,\exp(\beta T)\EE\,\exp(\beta \xi)\dd \suml_{i=0}^\infty \left(q^{i}\prod\limits_{j=1}^i\EE\,\left( \exp(\beta \{\xi_i|\overline{\EEE}_i\})\right)\right)\dd \frac{}{}
 \leqslant
\\
 \leqslant
 \dd \frac{\EE \,\exp(\beta t_1) \EE \,\exp(\beta t_1') \EE\,\exp(\beta \xi)}{1-qM(\beta\Theta)}=\mbox{\bf Exp}(\tau(b,b'),\beta).
\end{multline*}
\section{Conclusion}
After obtaining the bounds $\mbox{\bf Poly}(\tau(b,b'),\ell)$ or $\mbox{\bf Exp}(\tau(b,b'),\beta)$, we can integrate it by the stationary distribution $\tilde \PPP_B$ of the Markov process $B_t$.

This distribution is well-known (see \cite{Smith}): 
$$\tilde \PPP_B
\{[a,+\infty\}= \dd \frac{\displaystyle \int\limits _0^s (1-F(u))\,\mathrm{d} \, u} {\displaystyle \int\limits _0^\infty (1-F(u))\,\mathrm{d} \, u}= \displaystyle \dd \frac{\displaystyle \int\limits _0^s (1-F(u))\,\mathrm{d} \, u} {\EE\,\xi}. 
$$

For all $\ell\in(0,k-1]$ and all $\beta\in (0,\beta_0)$ these integrals converge to upper bounds in following inequalities
$$
\|\PPP(B_t)-\tilde\PPP_X\|_{TV} \leqslant \dd \frac{1}{t^{-\ell t}}\dd \intl_0^\infty \mbox{\bf Poly}(\tau(0,b'),\ell)\ud\tilde \PPP_B(b'),
$$

$$
\|\PPP(B_t)-\tilde\PPP_X\|_{TV} \leqslant e^{-\beta t}\dd \intl_0^\infty \mbox{\bf Exp}(\tau(0,b'),\beta)\ud \tilde\PPP_B(b'),
$$
accordingly.

These bounds are useful for the analysis and forecasting of the regenerative processes behaviour until the distribution of the process is close to the stationary invariant distribution. 

\end{exmp}

\end{document}